\documentclass[11pt, reqno, oneside]{article}
\usepackage{geometry} 
\usepackage[affil-it]{authblk}

\usepackage{amsmath, amsthm, amsfonts, amssymb, amsbsy, bigstrut, graphicx, enumerate,  upref, longtable, comment, booktabs, array, caption, subcaption}

\usepackage[breaklinks]{hyperref}

\usepackage[numbers, sort&compress]{natbib}

\newtheorem{thm}{Theorem}[section]

\theoremstyle{definition}

\newcommand{\pp}{\mathbb{P}}

\newcommand{\var}{\mathrm{Var}}

\numberwithin{equation}{section}

\usepackage[utf8]{inputenc}
\usepackage{amsmath}
\usepackage{amsfonts}
\usepackage{bbm}
\usepackage{mathtools}
\usepackage{tikz}
\usetikzlibrary{decorations.pathreplacing}
\usepackage{amsthm}
\usepackage[english]{babel}
\usepackage{fancyhdr}
\usepackage{amssymb}
\usepackage{enumitem}
\usepackage{qtree}
\usepackage{mathrsfs}

\newcommand{\E}{\mathbb{E}}

\renewcommand{\bar}{\overline}
\renewcommand{\P}{\mathbb{P}}

\begin{document}
\title{A Vershik--Kerov theorem for wreath products}
\author{Sourav Chatterjee\thanks{Departments of Mathematics and Statistics, Stanford University. \href{mailto:souravc@stanford.edu}{\tt souravc@stanford.edu}. 
}}
\affil{Stanford University}
\author{Persi Diaconis\thanks{Departments of Mathematics and Statistics, Stanford University. \href{mailto:diaconis@math.stanford.edu}{\tt diaconis@math.stanford.edu}. 
}}
\affil{Stanford University}


\maketitle


\begin{abstract}
Let $G_{n,k}$ be the group of permutations of $\{1,2,\ldots, kn\}$ that permutes the first $k$ symbols arbitrarily,  then the next $k$ symbols and so on through the last $k$ symbols. Finally the $n$ blocks of size $k$ are permuted in an arbitrary way. For $\sigma$ chosen uniformly in $G_{n,k}$, let $L_{n,k}$ be the length of the longest increasing subsequence in $\sigma$. For $k,n$ growing, we determine that the limiting mean of $L_{n,k}$ is asymptotic to $4\sqrt{nk}$. This is different from parallel variations of the Vershik--Kerov theorem for colored permutations.
\end{abstract}

\section{Introduction}
One of the most influential theorems of twentieth century probability is the Vershik--Kerov--Logan--Shepp solution of Ulam's problem: Determine the limiting mean of $L(\sigma)$, for $\sigma$ uniform in $S_n$, where $L(\sigma)$ is the length of the longest increasing subsequence in $\sigma$. This was proved independently by Vershik and Kerov~\cite{vershikkerov77, vershikkerov85} and Logan and Shepp~\cite{loganshepp77}, following substantial effort by Ulam~\cite{ulam61} and Hammersley~\cite{hammersley72}.  Later refinements by Baik, Deift and Johansson~\cite{baiketal99} determined the fluctuations and the limiting distribution
\[
\P\biggl(\frac{L(\sigma)-2\sqrt{n}}{n^{1/6}} \le x\biggr)\to TW_2(x)
\]
where $TW_2$ is the GUE Tracy--Widom distribution. This opened up a world of connections between probabilistic combinatorics and random matrix theory which is alive and thriving fifty years later.  A detailed history of these developments can be found in~\cite{aldousdiaconis95, romik15, aldousdiaconis99, corwin18}.

For reasons explained in Section \ref{sectwo}, we were interested in parallel results for the wreath product $G_{n,k} = S_k^n \rtimes S_n$. This acts as a subgroup of $S_{nk}$ by permuting symbols by an arbitrary $\gamma_1\in S_k$ for positions $\{1,\ldots,k\}$, an arbitrary $\gamma_2\in S_k$ for positions $\{k+1,\ldots,2k\}$, $\ldots,$ an arbitrary $\gamma_n\in S_k$ for positions $\{n(k-1)+1,\ldots,nk\}$. This is followed by permuting these $n$ blocks by an arbitrary $\eta\in S_n$. Let $\sigma = (\gamma_1,\ldots,\gamma_n;\eta)$. 

For example ($k=2$, $n=3$): $((1\, 2), (1)(2), (1\, 2); (3\, 1\, 2))$ permutes $1\, 2\, 3\, 4\, 5\, 6$ first to $2\, 1\, 3\, 4\, 6\, 5$ and then $6\, 5\, 2\,1\, 3\, 4$, with $L(\sigma)=3$.

The main result of this paper is: 

\begin{thm}\label{mainthm}
Let $n$ and $k$ be positive integers. Partition the numbers $1,2, \ldots,nk$ into $n$ blocks of size $k$ each. Randomly permute the blocks, and randomly permute within each block, all independently. Let $L_{n,k}$ be the length of the longest increasing subsequence in the resulting permutation. If $n, k\to \infty$ in such a way that $k/(\log n)^2 \to \infty$, then 
\[
\frac{L_{n,k}}{4\sqrt{nk}} \to 1 \ \text{ in probability.}
\]
Moreover, even if $k$ remains fixed or grows at any arbitrary rate as $n\to \infty$, we always have, for each $\epsilon >0$,
\[
\lim_{n\to \infty}\pp\biggl(\frac{L_{n,k}}{4\sqrt{nk}} \ge 1-\epsilon\biggr) =1.
\]
\end{thm}

Section \ref{sectwo} discusses background and other combinatorial probabilistic results for random elements of wreath products. The proof of Theorem \ref{mainthm} is in Section \ref{secthree}. We are unable to determine the growth rate of $L_{n,k}$ when $k$ is bounded or determine its fluctuations or limiting distributions. We wish we could ask our colleagues Vershik or Kerov.

\section{Combinatorial probability for wreath products}\label{sectwo}
\subsection{Background on wreath products}
Let $\Gamma \leqslant S_k$ and $H \leqslant S_n$ be subgroups. Their wreath product $\Gamma^n \rtimes H$ acts on $\{1,\ldots,nk\}$ as above. Familiar examples include:
\begin{itemize}
\item $B_n = S_2^n \rtimes S_n$ -- the hyperoctahedral group -- group of symmetries of an $n$-dimensional hypercube. This occurs as the group of centrally symmetric permutations in $S_{2n}$; these are the possible permutations with the two types of perfect shuffles~\cite{diaconisetal83}. In \cite{diaconisholmes98} it is shown how $B_n$ bijects with phylogenetic trees and the number of fixed points is used to bound rates of convergence of a walk natural in biology.
\item $C_k^n\rtimes S_n$ -- the generalized symmetric group. This may be seen as the group of $n\times n$ permutation matrices with the usual `ones' replaced by $k^{\mathrm{th}}$ roots of unity. It has its own Wikipedia page and of course $C_2^n\rtimes S_n = B_n$.
\item $S_k^n \rtimes S_n$ -- a maximal subgroup of $S_{kn}$ through the O'Nan--Scott theorem classifying primitive actions of the symmetric group~\cite{dixonmortimer96}. See \cite{diaconisetal08} for enumerative results and many further references.
\end{itemize}
More general wreath products occur throughout group theory (often as counter examples). It is natural to ask 
\begin{quote}
`Pick $\sigma \in \Gamma^n \rtimes H$ at random; what does it look like?'
\end{quote}
This is the familiar topic of combinatorial probability, a favorite topic of Anatoli Vershik and many collaborators.

A slew of first results are derived in \cite{diaconistung24}: The joint distribution of cycles, descents and inversions, and other features. Their motivation was the analysis of a new algorithm for generating random partitions of $n$. More generally, let $G$ be a finite group. The commuting graph of $G$ has vertex set $G$ an edge from $s$ to $t$ if $st=ts$. The natural Markov chain has transition matrix $K(s,t) = 1/|C_G(s)|$ (or $0$ if $s$ and $t$ do not commute). This is a reversible Markov chain on $G$ with stationary distribution $\pi(s)=1/|K(s)|$ where $K(s)$ is the conjugacy class containing $s$. Thus, running this chain and simply reporting the class gives a Markov chain with stationary distribution uniform on classes. When $G=S_n$, the classes are indexed by partitions of $n$, so the commuting graph walk gives a simple-to-run procedure for generating a random partition. Practical experiments show this may be the current `algorithm of choice'.

For $\sigma \in S_n$, $C_{S_n}(\sigma) = \prod_{i=1}^n C_i^{a_i}\rtimes S_{a_i}$, where $a_i$ is the number of cycles of length $i$ in $\sigma$. Proving things about rates of convergence requires understanding these wreath products. See \cite{diaconistung24, diaconishowes25} for details.

With this background, it is natural to ask about other features of $\sigma$. $L_n(\sigma)$ has generated such rich mathematics that this is a first further choice.

\subsection{Colored permutations}
There have been many variations and extensions of $L(\sigma)$ on $S_n$. One close to (but different from) the present study is the longest increasing subsequence of a `colored permutation'~\cite{tracywidom99, borodin99, rains98}. Here $\sigma$ is a permutation with each symbol decorated by one of $m$ colors (so there are $n!m^n$ such). $L(\sigma)$ is the length of the longest increasing subsequence in $\sigma$ with each symbol of the same color. The paper~\cite{borodin99} cited above proves analogues of the Baik--Deift--Johansson theorem for colored permutations. It develops deep connections to the distribution of eigenvalues of uniform elements of the unitary group and much else. We recommend \cite{borodin99} for background and developments.

When $n=2$, one talks about `signed permutations' which biject with the hyperoctahedral group $B_n$ above as permutations $\sigma$ of $\{-n, -(n-1),\ldots,-2,-1,1,2,\ldots, n\}$ which preserve central symmetry $\sigma(-x)=-\sigma(x)$. Here $-n,\ldots,n$ has its usual order with the natural definition of $L(\sigma)$. Tracy and Widom~\cite{tracywidom99} proved that:
\begin{thm}[\cite{tracywidom99}]
For $\sigma$ uniform in $B_n$,
\[
\P\biggl(\frac{L(\sigma) - 2\sqrt{2n}}{(2n)^{1/6}} \le x\biggr)\to  TW_2(2^{-2/3} x),
\] 
where $TW_2$ is the GUE Tracy--Widom distribution.
\end{thm}
Thus, for this representation of $B_n$, $L(\sigma)$ is about $2\sqrt{2n}$. Theorem \ref{mainthm} above does not work when $k=2$ but we conjecture that, in the wreath product action, $L(\sigma) \sim 3\sqrt{n}$.

The difference lies in there being two different actions of the same group(!). This is easiest to see through an example in $S_4$ with $S_2^2 \rtimes S_2$, a subgroup of order $8$.

\begin{table}[h]
\begin{center}
\caption{Wreath action}
\begin{tabular}{ccccccccc}
\toprule
$\sigma$ & $1234$ & $2134$ & $1243$ & $2143$ & $3412$ & $4312$ & $3421$ & $4321$\\
\midrule
$L(\sigma)$ & $4$ & $3$ & $3$ & $2$ & $2$ & $2$ & $2$ & $1$\\
\bottomrule
\end{tabular}
\end{center}
\end{table}

\begin{table}[h]
\begin{center}
\caption{Signed action}
\begin{tabular}{ccccccccc}
\toprule
$\sigma$ & $\bar{2}\bar{1}12$ & $2\bar{1}1\bar{2}$ & $21\bar{1}\bar{2}$ & $\bar{2}1\bar{1}2$ & $\bar{1}\bar{2}21$ & $\bar{1}2\bar{2}1$ & $12\bar{2}\bar{1}$ & $1\bar{2}2\bar{1}$\\
\midrule
$L(\sigma)$ & $4$ & $2$ & $1$ & $3$ & $2$ & $2$ & $2$ & $2$\\
\bottomrule
\end{tabular}
\end{center}
\end{table}

The distributions of $L(\sigma)$ are different. We have hopes that the natural wreath action might fit with nice mathematics and give a limit theory for $L(\sigma)$.

Perhaps it is worth explaining the problem of extending Theorem \ref{mainthm} to the case of bounded $k$. In the wreath product action, an increasing subsequence of the $n$ blocks of length $l$ gives an increasing subsequence of the full permutation of length at least $l$. We expect half of the size $2$ blocks to be increasing and so, using the Vershik--Kerov theorem, expect $L_{n,2}\sim \frac{3}{2}\cdot 2\sqrt{n} = 3\sqrt{n}$. However, there are {\it many} different increasing subsequences of length approximately $2\sqrt{n}$ in the $n$-block permutation, and there may be an anomalous one, just by chance, whose length is bigger than $3\sqrt{n}$ by a constant factor. When $k$ is a `bit large' (e.g., $k\gg (\log n)^2$ in Theorem \ref{mainthm}), we can prove sufficient concentration to rule out such anomalies.

\section{Proof of Theorem \ref{mainthm}}\label{secthree}
Fix some $n$ and $k$. For $i=1,\ldots,n$, let $N_i$ be the length of the longest increasing subsequence in block $i$. Let $N$ be the length of the longest increasing subsequence in the permutation of the blocks. Let $i_1<i_2<\cdots < i_N$ be the indices of the blocks in this longest increasing subsequence (and if there is  more than one longest increasing subsequence, choose one according to some predetermined rule). Then clearly, we can construct an increasing subsequence of numbers by putting together the longest subsequences in these blocks one after the other. Thus,
\begin{align}\label{lnkw}
L_{n,k} \ge W:= \sum_{j=1}^N N_{i_j}.
\end{align}
Note that $\E(N_i) = f(k)$ and $\var(N_i)=g(k)$, where $f(k)$ and $g(k)$ are the expected value and the variance of the length of the longest increasing subsequence in a uniform random permutation of $1,\ldots,k$. Moreover, $N_1,\ldots,N_n$ are independent, and $N$ is independent of $N_1,\ldots,N_n$. Thus,
\[
\E(W) = \E[\E(W|N)] = \E(N f(k)) = f(n) f(k),
\]
and
\begin{align*}
\var(W) &= \E(\var(W|N)) + \var(\E(W|N))\\
&= \E(N g(k)) + \var(N f(k))\\
&= f(n) g(k) + g(n) f(k)^2. 
\end{align*}
Now, we know that $f(n) \sim 2\sqrt{n}$ and $g(n) \sim Cn^{1/3}$ as $n\to \infty$, where $C$ is a positive constant whose value is known from~\cite[Theorem 1.2]{baiketal99}. Thus, $\E(W)\sim 4\sqrt{nk}$ and $\var(W) \sim 2C \sqrt{n} k^{1/3} + 4Ck n^{1/3}$ as $n,k\to \infty$. Now, if $n\to \infty$ with $k$ varying arbitrarily with $n$ (or remaining fixed), then $\sqrt{n} k^{1/3} = o(nk)$ and $k n^{1/3} = o(nk)$. Thus,
\[
\frac{\E(W)}{4\sqrt{nk}} \to 1, \ \ \ \var\biggl(\frac{W}{4\sqrt{nk}}\biggr) \to 0.
\]
Consequently, $W/(4\sqrt{nk}) \to 1$ in probability. Thus, by \eqref{lnkw}, we get that for any $\epsilon >0$,
\begin{align}\label{bd1}
\lim_{n,k\to \infty}\pp\biggl(\frac{L_{n,k}}{4\sqrt{nk}} \ge 1-\epsilon\biggr) =1.
\end{align}
(We remark here that the full force of \cite{baiketal99} is not really needed for the above argument. It suffices to observe that $g(n) = O(\sqrt{n})$, which follows, for example, from a simple argument via the Efron--Stein inequality.)
 
Next, fix some $\epsilon >0$ and some integers $n$ and $k$. Let $L_k$ denote the length of the longest increasing subsequence in a uniform random permutation of $1,\ldots,k$. By \cite[Theorem 7.1.2]{talagrand95}, we have that for any $u\ge 0$,
\[
\pp(L_k \ge M_k + u) \le 2\exp\biggl(-\frac{u^2}{4(M_k+u)}\biggr),
\]
where $M_k$ is any median of the law of $L_k$. Now, by Chebychev's inequality,
\[
\P(L_k \ge f(k) + 2\sqrt{g(k)}) \le \frac{g(k)}{4g(k)} = \frac{1}{4}. 
\]
This shows that we can find $M_k \le f(k) + 2\sqrt{g(k)} =: h(k)$. Consequently,
\begin{align*}
\pp(L_k \ge h(k) + u) &\le \pp(L_k \ge M_k + u) \\
&\le 2\exp\biggl(-\frac{u^2}{4(M_k+u)}\biggr)\\
&\le 2\exp\biggl(-\frac{u^2}{4(h(k)+u)}\biggr).
\end{align*}
Since $f(k)\sim 2\sqrt{k}$ and $g(k)\sim C k^{1/3}$ (or, with less effort, $g(k)=O(\sqrt{k})$), this shows that
\[
\P(L_k > 2\sqrt{k}(1+\epsilon)) \le 2 e^{-C_1\epsilon^2 \sqrt{k}},
\]
where $C_1$ is a universal constant. 
Now let $N$ and $N_i$ be as before. Let $E$ be the event that $N\le 2\sqrt{n}(1+\epsilon)$ and $N_i \le 2\sqrt{k}(1+\epsilon)$ for each $i$. Then by the above bound, 
\begin{align*}
\P(E^c) &\le \P(N > 2\sqrt{n}(1+\epsilon) + n\P(N_1> 2\sqrt{k}(1+\epsilon))\\
&\le 2 e^{-C_1\epsilon^2 \sqrt{n}} + 2 n e^{-C_1\epsilon^2 \sqrt{k}}.
\end{align*}
Note that $\P(E^c)\to 0$ if $n\to \infty$ and $k$ grows with $n$ at a rate faster than $(\log n)^2$. Now suppose that $E$ happens.  Let $i_1<i_2<\cdots < i_{L_{n,k}}$ be a longest increasing subsequence in the overall permutation of $1,\ldots,nk$. This must be the union of longest increasing subsequences in a sequence of the blocks, because otherwise we can increase its length by adding extra elements. Call this sequence of blocks $A$. Note that  $A$ must be an increasing subsequence of blocks (although $A$ may not be a {\it longest} increasing subsequence of blocks). Since $E$ has happened, $|A|\le 2\sqrt{n}(1+\epsilon)$ and $N_i\le 2\sqrt{k}(1+\epsilon)$ for each $i$. Thus,
\begin{align*}
L_{n,k} &= \sum_{i\in A} N_i \le |A| \max_{1\le i\le n} N_i \le 4\sqrt{nk}(1+\epsilon)^2.
\end{align*}
This proves that for any $\epsilon>0$,
\[
\lim_{n,k\to \infty, \, k\gg (\log n)^2} \P(L_{n,k} > 4(1+\epsilon)^2\sqrt{nk} ) = 0.
\]
Combined with \eqref{bd1}, this completes the proof.

\section*{Acknowledgements}
We thank Alexei Borodin for telling us about colored permutations and Nathan Tung for collaborative work which blends with the present project. We also thank the referee for a number of useful suggestions, and Lucas  Teyssier for an improvement of the proof of the theorem that relaxed the condition $k\gg (\log n)^4$ that we originally had, to the current condition $k\gg (\log n)^2$.

\bibliographystyle{abbrvnat}

\bibliography{myrefs}

\end{document}